\documentclass[a4paper, reqno, 12pt]{amsart}
\usepackage[english]{babel}
\usepackage[T1]{fontenc}
\usepackage{mathptmx,amsmath}
\usepackage{hyperref}
\usepackage{graphicx}

\usepackage{fullpage}

 \theoremstyle{remark}

\numberwithin{equation}{section}

\def\cen{\centerline}
\def\loc{\text{loc}}
\def\n{\noindent}
\def\al{\alpha}
\def\om{\omega}
\def\la{\lambda}
\def\va{\varphi}

\def\Om{\Omega}
\def\fr{\frac}

\def\psh{plurisubharmonic}
\def\ve{\varepsilon}
\def\de{\delta}
\def\PSH{\text{PSH}}
\def\nhd{neighborhood}

\def\De{\Delta}

\begin{document}
\title {A comparison principle for bounded plurisubharmonic functions on complex varieties in $\mathbb C^n$}
\author{Nguyen Quang Dieu and Sanphet Ounheuan}
\address{Department of Mathematics, Hanoi National University of Education,
136 Xuan Thuy street, Cau Giay, Hanoi, Vietnam}
\email{dieu$\_$vn@yahoo.com, sanphetMA@gmail.com}



\date{\today}

\maketitle
\section{Introduction}
\vskip0,3cm
\noindent
Let $D$ be a bounded domain in $\mathbb C^n.$ Denote by $\PSH (D)$ the cone of \psh\ functions on $D$ and 
$PSH(D) \cap L^\infty_{\loc} (D)$ (resp. $PSH(D) \cap L^\infty (D)$ 
the sub-cone of locally bounded (resp. bounded) \psh\ functions on $D$.
According to the fundamental work of Bedford and Taylor (see [BT1], [BT2], [BT3]), the complex Monge-Amp\`ere operator
$(dd^c)^n$ is well defined on $PSH(D) \cap L^\infty_{\loc} (D)$. This operator plays a prominent role in pluripotential theory just as the Laplace operator does in classical potential theory.
An important property of this operator is the following celebrated comparison principle due to Bedford and Taylor 
(see Theorem 4.1 in [BT1]).

\noindent
{\bf Theorem 1.1.} {\it Let $u, v \in PSH(D)\cap L^\infty (D)$ be such that 
$\varliminf\limits_{z \to \partial D} (u(z)-v(z)) \ge 0$. Then we have 
$$\int_{\{u<v\}} (dd^c v)^n \le \int_{\{u<v\}} (dd^c u)^n.$$}
\noindent
An analogous comparison principle was also obtained by Bedford (see Theorem 4.3 in [Be]) for bounded plurisubharmonic functions on open subsets of {\it complex spaces}. This result is the first inspiration for our
work. The other one comes from the following sharper form of Theorem 1.1 that was obtained
a few years later by Xing (see Lemma 1 in [Xi1]). 

\noindent
{\bf Theorem 1.2.} {\it Let $u, v \in PSH(D)\cap L^\infty (D)$ be such that 
$\varliminf\limits_{z \to \partial D} (u(z)-v(z)) \ge 0$. Then for any constant $r \ge 0$ and
$w_1, \cdots, w_n \in PSH(D)$ with $-1\le w_j< 0$ we have 
$$\begin{aligned}
&\frac1{(n!)^2} \int_{\{u<v\}} (v-u)^n dd^c w_1 \wedge \cdots\wedge dd^c w_n+
\int_{\{u<v\}} (r-w_1) (dd^c v)^n \\
&\le \int_{\{u<v\}} (r-w_1) (dd^c u)^n.
\end{aligned}$$}
Theorem 1.2 implies many important inequalities involving the complex Monge-Amp\`ere operator (see [Xi2] for details). Besides, this strong comparison principle
provides an effective tool in studying convergence problems for \psh\ functions and estimating capacity of small sets in pluripotential theory.
It should be noted, however, that in the extreme case $"r=\infty",$ Theorem 1.2 reduces to Theorem 1.1. Therefore, the essence of this version of the comparison principle lies at the other extreme $r=0.$
Along the development of energy classes for \psh\ functions (see [Ce]), there are variants of Theorem 1.2 that deal with 
\psh\ functions in Cegrell's classes, we could mention Theorem 4.7 in [KH] and Theorem 2 in [Xi2].

The aim of this note is twofold, first we generalize Theorem 1.2 to the context of bounded \psh\ functions on {\it complex varieties} in bounded domain of $\mathbb C^n$, and second we wish to relax somewhat the assumption on the boundary behavior of $u$ and $v.$ Another novelty of our work is to replace the expression $(v-u)^n$ in Theorem 1.2 by the composition of $v-u$ with a suitable real valued smooth function.

We now fix some notation and terminology that will be needed later on.
Given a connected complex variety $V$ of pure dimension $1 \le k \le n$ in 
a bounded domain $D$ in $\mathbb C^n,$ by $\PSH (V)$ (resp. $\PSH^{-} (V)$) we mean the set of \psh\ (resp. negative \psh) functions on $V.$ We defer to the next section for a brief account of plurisubharmonic functions on $V$ and the complex Monge-Amp\`ere operator on 
$\PSH (V) \cap L^{\infty}_{\loc} (V),$ the collection of locally bounded \psh\ functions on $V.$
A function $\chi: (0, \infty) \to (0, \infty)$ is said to be {\it $m-$increasing}, where $m \ge 1$ is an integer, 
if $\chi \in \mathcal C^m (0, \infty), \chi^{(j)}$ is increasing and non negative on $(0,\infty)$ for every 
$0 \le j \le m.$ For such a function $\chi$ and $0 \le j \le m,$ we set
\begin{equation} \label{chi}
\chi^{(j)}(0):= \lim_{t \to 0} \chi^{(j)} (t), P_m (\chi):= \sum_{j=0}^{m-1} \chi^{(j)} (0).
\end{equation}
A subset $E \subset \partial V$ is said to be {\it negligible} if there exists
$\psi \in PSH^{-} (V) \cap L^\infty_{\loc} (V)$ such that 
\begin{equation} \label{negli}
\underset{z \to \xi}{\lim} \psi(z)=-\infty, \ \forall \xi \in E.
\end{equation}
Our comparison principle reads as follows.
\vskip0,2cm
\noindent
{\bf Theorem 1.3.} {\it
Let $u,v \in PSH(V) \cap L^{\infty}_{\loc} (V)$ and $E \subset \partial V$ be a negligible set.
Assume that $u, v$ and $E$ satisfy the following conditions:

\noindent
(a) $\underset{z \in V}{\inf} (u(z)-v(z))>-\infty.$

\noindent
(b) $\varliminf\limits_{z \to \xi} (u(z)-v(z)) \ge 0$ for every $\xi \in (\partial V) \setminus E.$

Then for every integer $m$ with $1 \le m \le k$ and every $m-$increasing function $\chi: (0, \infty) \to (0, \infty)$ we have
$$\begin{aligned}
&\int_{\{u<v\}} \chi \circ (v-u)dd^c w_1 \wedge \cdots \wedge dd^c w_k +
\int_{\{u<v\}} (-w_1) \chi^{(m)} \circ (v-u) (dd^c v)^m \wedge dd^c w_{m+1} \wedge \cdots \wedge dd^c w_k\\
&\le \int_{\{u<v\}} (-w_1) \chi^{(m)} \circ (v-u) (dd^c u)^m \wedge dd^c w_{m+1} \wedge \cdots \wedge dd^c w_k
+P_m(\chi) \int_V dd^c w_1 \wedge \cdots \wedge dd^c w_k,
\end{aligned}$$
where $w_1, \cdots ,w_k \in PSH^{-} (V) \cap L^{\infty}_{\loc} (V)$ satisfying 
$w_j \ge -1$ for $2 \le j \le m$.}
\vskip0,3cm
\n
Let's point out that in the case where $V=D, m=k=n$ and $\chi (t)=t^n$, our comparison principle directly implies Theorem 1.2 even with a slightly better estimate, since $\chi^{(n)}\equiv n!<(n!)^2.$

The main ingredients in our poof are a smoothing method for \psh\ functions on complex varieties developed by Bedford in [Be] as well as integration by parts techniques demonstrated in [Xi1] and [KH]. 

The first application of our comparison principle is the following domination principle that was essentially due to Bedford and Taylor in the case where $V$ is an open domain in $\mathbb C^n$ and the exceptional set $E$ is empty (see Corollary 4.4 and Corollary 4.5 in [BT1]).

\noindent
{\bf Corollary 1.4.} {\it
Let $u,v \in PSH(V) \cap L^{\infty}_{\loc} (V)$ and $E \subset \partial V$ be as in Theorem 1.3.
Assume that for some $1 \le m \le k$ we have either $\underset{\{u<v\}}{\int} (dd^c u)^m \wedge \om^{k-m}=0$ or
$$(dd^c u)^m \wedge \om^{k-m} \le (dd^c v)^m \wedge \om^{k-m} \ \text{on the set}\ \{u<v\},$$
where $\om$ is the restriction of the K\"ahler form $dd^c \Vert z\Vert^2$ on $V.$
Then $u \ge v$ on $V.$} 
\vskip0,2cm
\noindent
The next consequence of Theorem 1.3 is a refinement of Theorem  4.3 in [Be].

\noindent
{\bf Corollary 1.5.} {\it Let $u,v \in PSH(V) \cap L^{\infty}_{\loc} (V)$ and $E \subset \partial V$ be as in Theorem 1.3. Then for every increasing continuous function $\chi: (0, \infty) \to (0,\infty)$ we have
$$\int_{\{u<v\}} \chi \circ (v-u)(dd^c v)^k \le \int_{\{u<v\}} \chi \circ (v-u) (dd^c u)^k.$$}
\noindent
We end up this section by presenting another consequence of Theorem 1.3 that offers a sufficient condition for convergence in capacity of a sequence in $PSH(V)\cap L^\infty (V).$ This result is similar in spirit to Theorem 3 in [Xi1] and Theorem 3.5 in [KH].
\vskip0,2cm
\noindent
{\bf Corollary 1.6.} {\it Let $u, \{u_j\} \subset PSH(V)\cap L^\infty (V).$ 
Let $\chi: (0, \infty) \to (0,\infty)$ be a increasing continuous function.
Assume that $u, u_j$ satisfy the following conditions:

\noindent
(a) $\underset{z \to \partial V}{\lim} (u(z)-u_j(z))=0$ for each $j \ge 1;$

\noindent
(b) $\underset{j \to \infty}{\lim} \underset{\{u_j<u\}}{\int} \chi \circ (u-u_j) d\vert \mu_j\vert=
\underset{j \to \infty}{\lim} \underset{\{u_j>u\}}{\int} \chi \circ (u_j-u) d\vert \mu_j\vert=0,$

where $\mu_j:=(dd^c u_j)^k-(dd^c u)^k.$

\n
Then $u_j \to u$ in capacity on $V$.}
 
\n
The conclusion of the above result says roughly that for each $\ve>0$ the capacities of the sets where the deviation of $u_j$ and $u$ is larger than $\ve$ tend to $0$ as $j \to \infty.$ 
Observe also that the sequence $\{u_j\}$ is not assumed to be locally {\it uniformly} bounded on $V.$

\vskip0,6cm
\noindent
{\bf Acknowledgments.} We are grateful to an anonymous referee for his(her) criticisms to an earlier version of this note. This work is written during a visit of the first named author to the Vietnam Institute for Advances in Mathematics (VIASM) in the winter of 2016. We would like to thank this institution for hospitality and financial support. Our work is also supported by the grant 101.02-2016.07 from the NAFOSTED program.
\section{Preliminaries}
We first recall elements of pluripotential theory on complex varieties in $\mathbb C^n.$ The main focus is the complex Monge-Amp\`ere operator and its continuity property. For more details we refer the reader to [Be].
 
Let $V$ be a connected complex variety of pure dimension $k$ in a bounded domain 
$D \subset \mathbb C^ n (n  \ge 2, 1 \le k \le n).$ Thus, $V$ is locally the common zero sets of holomorphic functions on open subsets of $D.$
We denote by $V_r$ the set of regular points of $V.$ Hence $V_r$ is the largest (possibly disconnected) complex manifold of dimension $k$ included in $V.$ The singular locus of $V$ is then denoted by $V_s:=V \setminus V_r.$
A function $u: V \to [-\infty, \infty)$ is said to be \psh\ if $u$ is locally the restriction (on $V$) of \psh\ functions on an open subset of $D$. Notice that we regard the function identically $-\infty$ as plurisubharmonic.
A fundamental result of Fornaess and Narasimhan (see Theorem 5.3.1 in [FN]) asserts that an upper semicontinuous function $u: V \to [-\infty, \infty)$ is \psh\  if and only if the restriction of $u$ on every analytic curve in $V$ is subharmonic.
This powerful result implies immediately the nontrivial facts that  plurisubharmonicity is preserved under pointwise decreasing convergence. We write $PSH (V)$ for the cone of \psh\ functions on $V$ and $PSH^-(V)$ denotes  the sub-cone of negative \psh\ functions on $V$. 
According to a fundamental result of Lelong (see p.32 in [GH]), the set $V_r$ has {\it finite} volume near every point of $V_s.$ Therefore each $u \in PSH(V) \cap L^{\infty}_{\text{loc}} (V)$ defines a current of bidegree $(0,0)$ on $V.$
Hence we may regard $dd^c u$ as a {\it positive} closed current of bidegree $(1,1)$ on $V.$  

Next, we turn to the complex Monge-Amp\`ere operator for locally bounded \psh\ functions on $V.$
According to Bedford in [Be], the complex Monge-Amp\`ere operator
$$(dd^c )^k: PSH(V) \cap L^{\infty}_{\text{loc}} (V) \to M_{k,k} (V),$$
where $M_{k,k} (V)$ denotes the collection of Radon measures on $V,$ may be defined in the usual way on the regular locus $V_r$ of $V$ as in [BT1]. Namely, given $u \in PSH(V) \cap L^{\infty}_{\text{loc}} (V),$ we define inductively on $V_r$ the following currents
$$(dd^c u)^m:=dd^c (u(dd^c u)^{m-1}), 1 \le m \le k,$$
and the measure $(dd^c u)^k$ extends "by zero" through the singular locus $V_s,$ i.e., for Borel sets
$E \subset V$
\begin{equation}
\label{monge}
\int_E (dd^c u)^k: =\int_{E \cap V_r} (dd^c u)^k.
\end{equation}
\noindent
For a Borel subset $E$ of an open set $\Om \subset V,$ the capacity of $E$ relative to $\Om$ is defined by
$$C(E,\Om)=\sup\big\{\int_E (dd^c u)^k: u \in \PSH(\Om): -1 \le u<0\big\}.$$
The above definition makes sense since $\psi (z):=\fr{\Vert z\Vert^2}M-1 \in PSH(V)$ and satisfies $-1\le \psi<0,$
where $M:=\sup \{\Vert z\Vert^2: z \in D\}<\infty.$

Obviously, by (\ref{monge}), the singular locus $V_s$ has zero capacity, i.e., $C(V_s \cap \Om, \Om)=0$ for every open subset $\Om$ of $V.$
The following basic result of Bedford (Lemma 3.1 in [Be]) asserts that $V_s$ actually has {\it outer} capacity zero.

\noindent
{\bf Lemma 2.1.} {\it For every open subset $\Om$ of $V$ and
every $\ve>0$, there exists an open \nhd\ $U$ of $V_s$ in $\Om$ such that
$C(U, \Om)<\ve.$}

\noindent
From Lemma 2.1, we see immediately that (\ref{monge}) in fact defines $(dd^c u)^k$ as a Radon measure on $V$
for each  $u \in PSH(V) \cap L^{\infty}_{\text{loc}} (V).$
Moreover, the following version of the Chern-Levine-Nirenberg inequality holds true on $V:$ For every relatively compact open 
subset $V' \subset V$ and every Borel subset $E$ of $V'$ we have
\begin{equation} \label{CLN}
\int_E (dd^c u)^k \le \la \Vert u\Vert_{V'}^k,
\end{equation}
where $\la \ge 0$ is a finite constant depends only on $E, V'.$

More generally, if $u_1, \cdots, u_k \in PSH(V) \cap L^\infty_{\loc} (V)$ then by local polarization in the symmetric linear form $dd^c u_1 \wedge \cdots \wedge dd^c u_k,$ e.g, for every $U$ open relatively 
compact subset of $V$ we set
$$2dd^c u'_1 \wedge dd^c u'_2=dd^c (u'_1+u'_2)^2-dd^c ({u'_1}^2)- dd^c ({u'_2}^2),$$
where $u'_1=u_1-\inf_U u_1, u'_2:=u_2-\inf_U u_2,$
we see that $dd^c u_1 \wedge \cdots \wedge dd^c u_k$ defines a Radon measure on $V$ as well.
Another point to stress is that, all the local analysis in the fundamental work [BT1] carries over $V_s$. For instance, from Theorem 3.5 in [BT1] and Lemma 2.1 we conclude that every $u \in \PSH(V)$ is quasi-continuous on $V,$ i.e., for every $\ve>0$ 
there exists an open subset $V_\ve$ of $V$ such that $C(V_\ve, V)<\ve$ and $u$ is continuous on $V\setminus V_\ve$.
It follows, using Dini's lemma, that every sequence $\{u_j\} \in \PSH(V)$ that converges monotonically to $u \in \PSH(V)$ must converge
locally quasi-uniformly i.e., given a relatively compact open subset $V'$ of $V$ and $\ve>0$, there exists an open subset
$V'_\ve \subset V'$ such that $C(V'_\ve, V)<\ve$ and $u_j$ converges uniformly to $u$ on $V'\setminus V'_\ve.$ 

We claim no originality for the following result about convergence of certain currents on $V$.

\noindent
{\bf Proposition 2.2.} {\it Let $p,q,r$ be non-negative integers and
$\{u_{1,j}\}, \cdots, \{u_{p,j}\},
\{v_{1,j}\}, \cdots, \{v_{q,j}\},$ 

\noindent
$\{w_{1,j}\}, \cdots, \{w_{r,j}\}$ be sequences 
in $\PSH(V)$ that decrease pointwise to $u_1, \cdots, u_p,
v_1, \cdots, v_q,$ 

\noindent
$w_1, \cdots, w_r \in PSH(V) \cap L^\infty_{loc} (V).$ For each $j,$ define the current
$$T_j:=du_{1,j}\wedge \cdots \wedge du_{p,j}\wedge d^cv_{1,j}\wedge \cdots\wedge d^cv_{q,j}\wedge dd^c w_{1,j}\wedge\cdots\wedge dd^c w_{r,j}.$$
Then the following assertions hold true:

\noindent
(a) The total variation of norms $\vert T_j\vert$ of $T_j$ are uniformly bounded on compact sets of $V;$

\noindent
(b) $T_j$ converges weakly to 
$T:=du_1\wedge \cdots \wedge du_p\wedge d^cv_1\wedge \cdots\wedge d^cv_q\wedge dd^c w_1\wedge\cdots\wedge dd^c w_r;$

\noindent
(c) If $\{\psi_j\}, \psi$ are quasi-continuous functions on $V$ which are locally uniformly bounded 
and if $\psi_j$ converges locally quasi-uniformly to $\psi$ then 
$\psi_jT_j$ converges weakly to $\psi T.$}
\vskip0,2cm
\noindent
{\it Proof.} Given a relatively compact open subset $V'$ of $V$ and a compact subset $K$ of $V' \cap V_r,$ by {\it the proof} of Lemma 2.2 of [BT3] where a stronger version of the Chern-Levine-Nirenberg inequality (\ref{CLN}) is established, we can find a constant $\gamma>0$ that depends only on the sup norms on $V'$ of $u_1, \cdots, u_p,
v_1, \cdots, v_q, w_1, \cdots, w_r$
such that for $j$ large enough we have
\begin{equation} \label{CLNN}
\vert T_j \vert (K) \le \gamma C(K,V').
\end{equation}
Using Lemma 2.1 we see that the above inequality holds true for {\it any} compact subset $K$ of $V'.$
This proves the statement (a).
Next, (\ref{CLNN}) also implies that $\vert T_j \vert$ put {\it uniformly} small mass on sets having small capacity, i.e., given a relatively compact open subset $V'$ of $V$ and $\ve>0$, there exists $\de>0, j_0 \ge 1$ such that for every compact subset
$K \subset V'$ with $C(K, V')<\de$ we have 
$\vert T_j\vert (K)<\ve$ for $j \ge j_0.$ This fact together with Proposition 2.3 in [BT3] implies (b).
Finally (c) follows from an easy adaptation of the proof $(4) \Rightarrow (1)$ in Theorem 3.2 of [BT2] (see also Theorem 2.6 in [BT3]).
\vskip0,2cm
Now we discuss the problem of smoothing \psh\ functions on $V$.
In the case where $V$ is {\it Stein,} i.e., there exists a strictly \psh\ exhaustion for $V$, we can approximate every element $u \in \PSH(V)$ from above by a decreasing sequence of $\mathcal C^\infty$ smooth strictly \psh\ functions on $V$ (see Theorem 5.5 in [FN]).
For a general $V$, such a smoothing may not be possible even on domains in $\mathbb C^n$ 
(see p. 297 in [Be] for a counterexample of Fornaess).
So we have to be content with the following smoothing method devised by Bedford (see p. 299 in [Be]). Namely, let $\psi \in PSH(V)$ be given, and let $\mathcal U:=\{U_l\}$ be an open covering of $V$ such that 
for each $l$ there is an open subset $\tilde U_l$ of $D, U_l$ is a complex variety of $\tilde U_l,$ and there exists
$\tilde \psi_l \in \PSH(\tilde U_l)$ with $\tilde \psi_l=\psi$ on $U_l.$
Next, we let $\{\chi_l\}$ be a partition of unity subordinate to $\tilde {\mathcal U}:= \{\tilde U_l\}$. For each $l,$ after  taking convolution 
$\tilde \psi_l$ with standard radial smoothing kernels $\rho_\de$ on $\mathbb C^n$, we obtain a smoothing $\tilde \psi_{l,\de}$ which is smooth and \psh\ on a \nhd\ of $\text{supp}\ \chi_l$ for $\de>0$ small enough.
Now our smoothing is obtained as the sum
\begin{equation} \label{smoothing}
\psi^\de:= \sum \chi_l \tilde \psi_{l,\de}.
\end{equation}
It is clear that $\psi^\de$ is smooth on a \nhd\ of $K$ for every compact subset $K$ of $V$.
Moreover, $\psi^\de \downarrow \psi$ on $V$ as $\de \to 0.$
In general, $\psi^\de \not\in PSH(V)$. However, as we will see below, these smoothings are nice enough to make continuity of the complex Monge-Amp\`ere operator possible.
More precisely, let $u_1, \cdots, u_k \in PSH(V)\cap L^\infty_{\loc} (V).$  Choose a common covering 
$\mathcal U:= \{U_l\}$ of $V$ and a partition of unity $\{\chi_l\}$ subordinate to $\tilde {\mathcal U}$ for all \psh\ functions $u_1, \cdots, u_k$. Then we have the following approximation result which is implicitly contained in [Be].
\vskip0,2cm
\noindent
{\bf Proposition 2.3.} {\it Let $\{f_j\}, f$ be
locally uniformly bounded, quasi-continuous functions on $V$. Assume that $\{f_j\}$ converges locally 
quasi-uniformly to $f.$
Then for every sequence $\{\de_j\}\downarrow 0,$ 
the currents $f_jdd^c u_1^{\de_j} \wedge \cdots\wedge dd^c u_k^{\de_j}$ converges weakly to $fdd^c u_1 \wedge \cdots \wedge dd^c u_k$ as $j \to \infty.$
}
\vskip0,2cm
\noindent
{\it Proof.} For each $p, 1 \le p \le k$ and $j \ge 1,$ by (\ref{smoothing}) we have the following equalities on $V$
$$u_p^{\de_j}= \sum_l \chi_l.{(\tilde u_p)}_{l, \de_j}, u_p= \sum_l \chi_l.{(\tilde u_p)}_l.$$
For simplicity of notation,  put $v_{p,l,j}:= {(\tilde u_p)}_{l, \de_j}$ and 
$v_{p,l}:={(\tilde u_p)}_l.$
Then $v_{p,l,j} \downarrow v_{p,l}$ on $U_l$ as $j \to \infty.$
By direct computation, we expand 
$$dd^c u_p^{\de_j}= \sum_l v_{p,l,j}dd^c \chi_l+d\chi_l \wedge d^c v_{p,l,j}+ 
dv_{p,l,j} \wedge d^c \chi_l+\chi_l dd^c v_{p,l,j},$$
$$dd^c u_p= \sum_l v_{p,l}dd^c \chi_l+d\chi_l \wedge d^c v_{p,l}+ 
dv_{p,l} \wedge d^c \chi_l+\chi_l dd^c v_{p,l}.$$
Hence each of the currents $dd^c u_1^{\de_j} \wedge \cdots\wedge dd^c u_k^{\de_j}$ and $dd^c u_1 \wedge \cdots \wedge dd^c u_k$ is 
the sum of $4^k$ smooth forms of bidegree $(k,k).$ Moreover, each of them, by an abuse of notation, can be represented as
$$\va_I w_{q_1,j}\cdots w_{q_r,j} dw_{q_{r+1},j}\wedge \cdots dw_{q_s,j} 
\wedge d^c w_{q_{s+1},j} \wedge \cdots\wedge d^c w_{q_t,j}\wedge 
dd^c w_{q_{t+1},j}\wedge \cdots \wedge dd^c w_{q_\al,j},$$
and
$$\va_I w_{q_1}\cdots w_{q_r} dw_{q_{r+1}}\wedge \cdots dw_{q_s} 
\wedge d^c w_{q_{s+1}} \wedge \cdots\wedge d^c w_{q_t}\wedge 
dd^c w_{q_{t+1}}\wedge \cdots \wedge dd^c w_{q_\al},$$
respectively, where 
$$\{w_{q_1, j}, \cdots, w_{q_\al,j}\} \subset \{v_{p,l,j}: 1 \le p \le k\},
\{w_{q_1}, \cdots, w_{q_\al}\} \subset \{v_{p,l}: 1 \le p \le k\},
I=(q_1,\cdots,q_\al)$$ 
and $\va_I$ are smooth forms with compact support that involves only on $\chi_l.$

By the assumptions on $\{f_j\}$ and $f$, we may use Proposition 2.2 (c) to conclude that 
for each $I=(q_1, \cdots, q_\al)$ the sequence of currents 
$$f_j \va_I w_{q_1,j}\cdots w_{q_r,j} dw_{q_{r+1},j}\wedge \cdots dw_{q_s,j} 
\wedge d^c w_{q_{s+1},j} \wedge \cdots\wedge d^c w_{q_t,j}\wedge 
dd^c w_{q_{t+1},j}\wedge \cdots \wedge dd^c w_{q_\al,j}$$
converges weakly to the current
$$f\va_I w_{q_1}\cdots w_{q_r} dw_{q_{r+1}}\wedge \cdots dw_{q_s} 
\wedge d^c w_{q_{s+1}} \wedge \cdots\wedge d^c w_{q_t}\wedge 
dd^c w_{q_{t+1}}\wedge \cdots \wedge dd^c w_{q_\al},$$
as $j \to \infty.$ By taking the sum over $I$ we obtain the desired conclusion.
\vskip0,4cm
\noindent
Our final auxiliary result concerns approximation of $m-$increasing functions by smooth ones.
\vskip0,2cm
\noindent
{\bf Lemma 2.4.} {\it Let $m \ge 1$ be an integer and $\chi: (0, \infty) \to (0, \infty)$ be a $m-$increasing function. 
Then there exists a sequence $\{\chi_j\}$ of $m-$increasing
$\mathcal C^\infty-$smooth functions  such that $\{\chi_j^{(l)}\}$ converges locally uniformly on $[0, \infty)$ to $\chi^{(l)}$ for $0 \le l \le m.$}
\vskip0,2cm
\noindent
{\it Proof.} Set $\va(t)=\chi^{(m)} (t)$ for $t \ge 0$ and $\va(t):=\chi^{(m)}(0)$ for $t<0.$ Then $\va$ is real valued, continuous and increasing on $\mathbb R.$ By taking convolution of $\va$ with approximate of identity, we obtain a sequence $\{\va_j\}$ of $\mathcal C^\infty-$smooth increasing functions on $\mathbb R$ that converge locally uniformly to $\va$.
Now for each $j,$ we define inductively on $(0, \infty)$ the following functions
$$\va_{j,0}:=\va_j, \va_{j,l}(t):= \int_0^t \va_{j,l-1}(x)dx+\chi^{(m-l)} (0), 1 \le l \le m.$$
Then we have $\va'_{j,l}=\va_{j,l-1}$ for $1 \le l \le m.$ 
Hence 
$\va_{j,m}^{(l)}=\va_{j,m-l}$ for $1 \le l \le m.$ 
Moreover, we can show by induction that 
$\{\va_{j,m}^{(l)}\}$ converges locally uniformly on $[0,\infty)$ to 
$\chi^{(l)}$ as $j \to \infty$ for $0 \le l \le m.$ It follows that $\chi_j:=\va_{j,m}$ is the sequence we are searching for.
\vskip0,4cm
\noindent
\section{Strong comparison principle}
\vskip0,2cm
\noindent
We start with the following simple facts that will be needed in examining certain integration by parts formulas. 
In the next two lemmas, we will denote by $\va$ a real valued $\mathcal C^2-$ smooth function defined on $(0, \infty).$

\noindent
{\bf Lemma 3.1.} {\it Let $u, v \in PSH(V) \cap L^\infty_{loc} (V)$ with $u<v$ on $V$. 
Then the following assertions hold in the sense of currents on $V:$

\noindent
(a) $dd^c (\va \circ (v-u))=\va' \circ (v-u) dd^c(v-u)+\va''(v-u)d(v-u)\wedge d^c (v-u).$ 

\noindent
(b) If $\va'' \ge 0$ on $(0, \infty)$ then $dd^c (\va \circ (v-u)) \ge \va' \circ (v-u)dd^c (v-u).$
}

\noindent
{\it Proof.} (a) Fix $a \in V$, it suffices to show the above identity in a small \nhd\ of $a$ in $V.$
Then we can find a ball $\mathbb B \subset D$ around $a$ and \psh\ functions $\tilde u, \tilde v$ on $\mathbb B$
such that $\tilde u|_{\mathbb B \cap V}=u, \tilde v|_{\mathbb B \cap V}=v$.
By considering $\max \{\tilde u, \tilde v\}$ instead of $\tilde v$, we can assume $\tilde u < \tilde v$ on $\mathbb B$.
By taking convolutions of $\tilde u$ and $\tilde v$ with standard radial smoothing kernels $\rho_\de$ on $\mathbb C^n$
and shrinking $\mathbb B,$ we obtain smooth \psh\ functions 
$u_\de, v_\de$ on $\mathbb {B}$ such that
$u_\de<v_\de$ and $u_\de \downarrow u, v_\delta \downarrow v$ on $\mathbb B \cap V$ as $\delta \downarrow 0.$
By direct computation we obtain for each $\de>0$ 
$$dd^c (\va \circ (v_\de-u_\de))=\va' \circ (v_\de-u_\de) dd^c(v_\de-u_\de)+\va'' \circ (v_\de-u_\de)d(v_\de-u_\de)\wedge d^c (v_\de-u_\de).$$
Since $u_\de, v_\de$ are uniformly bounded on $\mathbb B$, by letting $\de \downarrow 0$ and applying Lebesgue dominated convergence theorem we obtain the desired equality.

\noindent
(b) The inequality then follows directly from (a) and the fact that $d(v-u)\wedge d^c (v-u)$ is a non-negative $(1,1)$ current. 
\vskip0,2cm
\noindent
The following integration by parts formula plays a crucial role in the proof of Lemma 3.3.
Its proof requires all the machinery developed in the preceding section.

\noindent
{\bf Lemma 3.2.} {\it Let $u, v, w_1, \cdots, w_k \in PSH(V)\cap L^\infty_{\loc} (V).$
Assume that $u \le v$ on $V, u=v$ outside a compact subset $K$ of $V.$
Then for all real number $\ve>0$ and open sets $V'$ such that $K \subset V' \subset\subset V$ 
we have
$$\begin{aligned}
\int_{V'} \va \circ (v+\ve-u)dd^c w_1\wedge \cdots \wedge dd^c w_k&=\int_{V'} w_1 dd^c (\va \circ (v+\ve-u)) \wedge dd^c w_2\wedge \cdots \wedge dd^c w_k\\
&+\va(\ve)\int_{V'} dd^c w_1 \wedge \cdots \wedge dd^c w_k.
\end{aligned}$$
}
\noindent
{\it Proof.} For $\de>0$ small enough, following (\ref{smoothing}),
we let $u^\de, v^\de, w_1^\de, \cdots, w_k^\de$ be smoothing of 
$u, v, w_1, \cdots, w_k,$ respectively. Notice that the covering $\{U_j\}$ and the partition of unity $\{\chi_j\}$ can be chosen to be common for all these \psh\ functions. By the assumption we have $u^\de=v^\de$ on $V' \setminus K.$
In addition, as in the proof of Lemma 3.1, we may arrange so that 
$u^\de \le v^\de$ on $V'$ for every $\de.$ 
To simplify notation, we set 
$$T^\de:=dd^c w_2^\de \wedge \cdots \wedge dd^c w_k^\de, T:=dd^c w_2 \wedge \cdots \wedge dd^c w_k,
\va^\de: =\va \circ (v^\de+\ve-u^\de)-\va(\ve).$$
Notice that $\va^\de=0$ on a small \nhd\ of $\partial V'$.
Hence, an application of Stoke's theorem for smooth forms
on the complex variety $V'$ (see p.33 in [GH]) yields
$$\int_{V'} \big [\va^\de dd^c w_1^\de -w_1^\de dd^c \va^\de \big] \wedge T^\de 
=\int_{V'} d\big [\va^\de d^c w_1^\de \wedge T^\de -w_1^\de d^c \va^\de \wedge T^\de\big]=0.$$
It follows that 
\begin{equation} \label{Stoke}
\int_{V'} \big [\va \circ (v^\de+\ve-u^\de)-\va(\ve)\big ] dd^c w_1^\de \wedge T^\de
=\int_{V'} w_1^\de dd^c (\va \circ (v^\de+\ve-u^\de)) \wedge T^\de.
\end{equation}
We now consider the limits of both sides of (\ref{Stoke}) as $\de \downarrow 0.$ For the left hand side, set 
$$\mu^\de:= [\va \circ (v^\de+\ve-u^\de)-\va(\ve)] dd^c w_1^\de\wedge T^\de, 
\mu:= [\va \circ (v+\ve-u)-\va(\ve)]dd^c w_1 \wedge T.$$
Then $\mu^\de$ and $\mu$ are real measures on $V'$ that vanish outside $K$.
Observe that the functions $\va\circ (v^\de+\ve-u^\de)$ are continuous on $V,$ locally uniformly bounded and converges to 
$\va \circ (v+\ve-u)$ locally quasi-uniformly on $V'$. Hence, by Proposition 2.3  we deduce that $\mu^\de$ converge weakly to $\mu$ as $\de \downarrow 0.$ We claim that $\mu^\de (V') \to \mu(V')$ as $\de \downarrow 0.$ Indeed, fix $\ve>0.$
By the Chern-Levine-Nirenberg inequality (\ref{CLN}) we can choose an open subset $V''$ of $V$ such that
$K \subset V'' \subset \subset V$ and 
$\vert \mu \vert (V'' \setminus K)<\ve, \vert \mu^\de\vert (V'' \setminus K)<\ve$ for $0<\de<\de_0$ small enough.
Let $f$ be a continuous function on $V"$ with compact support such that 
$0 \le f \le 1, f=1$ on $K$.
Then we have
$$\vert \mu^\de (V')- \mu (V')\vert =\vert \mu^\de (K)-\mu (K)\vert 
\le \big \vert \int_{V''} f d\mu^\de-\int_{V''} f d\mu \big \vert+2\ve.$$
It follows that 
$$\varlimsup\limits_{\de \to 0} \vert  \mu^\de (V')- \mu (V')\vert \le 2\ve.$$
This proves our claim since $\ve>0$ can be chosen to be arbitrarily small.
Hence 
\begin{equation} \label{Stoke1}
\lim_{\de \to 0} \int_{V'} [\va \circ (v^\de+\ve-u^\de)-\va(\ve)] dd^c w_1^\de\wedge T^\de =
\int_{V'} [\va \circ (v+\ve-u)-\va(\ve)] dd^c w_1\wedge T.
\end{equation}
Similarly, for the right hand side of (\ref{Stoke}) we define the following currents on $V'$
$$\begin{aligned}
\mu'^{\de}:&= w_1^\de [\va' \circ (v^\de+\ve-u^\de)dd^c(v^\de-u^\de)+\va''\circ(v^\de+\ve-u^\de)d(v^\de-u^\de)\wedge d^c (v_\de-u_\de)\big ] \wedge T^\de,\\
\mu':&= w_1 [\va'\circ (v+\ve-u)dd^c(v-u)+\va''\circ (v+\ve-u)d(v-u)\wedge d^c (v-u)\big] \wedge T.
\end{aligned}$$
By repeating the same reasoning as above we have $\mu'^\de (V)\to \mu' (V)$ as $\de\downarrow 0.$
Therefore, by applying Lemma 3.1 (a) we obtain
\begin{equation} \label{Stoke2}
\lim_{\de \to 0} \int_{V'} w_1^\de dd^c (\va \circ (v^\de+\ve-u^\de)) \wedge T^\de=
\int_{V'} w_1 dd^c (\va \circ (v+\ve-u)) \wedge T.
\end{equation}
Combining (\ref{Stoke}), (\ref{Stoke1}) and (\ref{Stoke2}) together we obtain
$$\int_{V'} \big [\va \circ (v+\ve-u)-\va(\ve)\big ] dd^c w_1 \wedge \cdots \wedge dd^c w_k
=\int_{V'} w_1 dd^c (\va \circ (v+\ve-u)) \wedge dd^c w_2 \wedge \cdots \wedge dd^c w_k.$$
Finally, by (\ref{CLN}) we have $0 \le \int_{V'} dd^c w_1 \wedge \cdots \wedge dd^c w_k <\infty$, so after
rearranging the above equation we obtain the desired conclusion.
\vskip0,4cm
\noindent 
The next lemma, a special case of Theorem 1.3, is the key step in our proof. It is somewhat inspired from Lemma 3.3 in [KH].
\vskip0,2cm
\noindent
{\bf Lemma 3.3.} {\it Let $u, v\in \PSH (V) \cap L^\infty (V)$ be such that $u\le v$ on $V$ and 
$u=v$ outside a compact subset $K$ of $V$.
Then for any integer $1 \le m \le k$ we have
$$\begin{aligned}
&\int_V \chi \circ (v-u)dd^c w_1 \wedge \cdots \wedge dd^c w_k +
\int_V (-w_1) \chi^{(m)} \circ (v-u)(dd^c v)^m dd^c w_{m+1} \wedge \cdots \wedge dd^c w_k\\
&\le \int_V (-w_1) \chi^{(m)} \circ (v-u)(dd^c u)^m dd^c w_{m+1} \wedge \cdots \wedge dd^c w_k
+P_m(\chi) \int_V dd^c w_1 \wedge \cdots \wedge dd^c w_k,
\end{aligned}$$
where $w_1, \cdots ,w_k \in PSH (V) \cap L^\infty_{\loc} (V)$ satisfying $w_j<0$ for $1 \le j \le m$ and
$w_j \ge -1$ for $2 \le j \le m$.}
\vskip0,2cm
\noindent
{\it Proof.} For the ease of notation, we set 
$$T:= dd^c w_1 \wedge \cdots \wedge dd^c w_m, T':=dd^c w_{m+1} \wedge \cdots \wedge dd^c w_k.$$
Let $V'$ be a relatively compact connected open subset of $V$ such that $K \subset V'$. It follows that $u=v$ on a small
\nhd\ of $\partial V'.$ 
We are now aiming at the following estimate
$$\int_{V'} \chi \circ (v-u) T \wedge T' +
\int_{V'} (-w_1) \chi^{(m)} \circ (v-u)(dd^c v)^m \wedge T' $$
\begin{equation}\label{estimate}
\le \int_{V'} (-w_1) \chi^{(m)} \circ (v-u)(dd^c u)^m \wedge T'+P_m (\chi) \int_{V'} T \wedge T'.
\end{equation}
For this, we first assume that $\chi \in \mathcal C^{m+1} (0,\infty).$
Then $\chi^{(j)} \ge 0$ on $(0, \infty)$ for each $0 \le j \le m+1.$
Now by using the integration by parts formula (Lemma 3.2) and Lemma 3.1(b) (while noting that $w_m<0$) 
we obtain
$$\begin{aligned}
&\int_{V'} \chi \circ (v+\ve-u) T \wedge T' \\
&=\int_{V'} w_m dd^c (\chi \circ (v+\ve-u)) dd^c w_1 \wedge \cdots \wedge dd^c w_{m-1} \wedge T' +\chi(\ve)\int_{V'} T \wedge T'\\
& \le \int_{V'} w_m \chi' \circ (v+\ve-u) dd^c (v-u) \wedge dd^c w_1 \wedge \cdots \wedge dd^c w_{m-1} \wedge T'+
\chi(\ve)\int_{V'} T \wedge T'\\
& \le \int_{V'} \chi' \circ (v+\ve-u)dd^c u \wedge dd^c w_1 \wedge \cdots \wedge dd^c w_{m-1} \wedge T'
+\chi(\ve)\int_{V'} T \wedge T'.
\end{aligned}
$$
Here the last inequality follows from the fact that $w_m dd^c (v-u) \le dd^c u,$ which in turns 
is a consequence of the assumption that $-1 \le w_m<0$ on $V.$
Continuing this process $(m-2)$ more times we get
$$
\int_{V'} \chi \circ (v+\ve-u) T \wedge T' \le 
\int_{V'} \chi^{(m-1)} \circ (v+\ve-u)(dd^c u)^{m-1} \wedge dd^c w_1 \wedge T'
+(\sum_{j=0}^{m-2} \chi^{(j)} (\ve))\int_{V'} T \wedge T'.$$
Next, since $\chi \in \mathcal C^{m+1} (0,\infty)$ we may apply Lemma 3.2 and Lemma 3.1(b) again (while noting that $w_1<0$ on $V$) to get
$$\begin{aligned} 
\int_{V'} \chi^{(m-1)} \circ (v+\ve-u)(dd^c u)^{m-1} \wedge dd^c w_1 \wedge T' 
&=\int_{V'} w_1 dd^c (\chi^{(m-1)} \circ (v+\ve-u))\wedge (dd^c u)^{m-1} \wedge T' \\
&+\chi^{(m-1)}(\ve) \int_{V'} dd^c w_1 \wedge \cdots \wedge dd^c w_k \\
&\le \int_{V'} w_1 \chi^{(m)} \circ (v+\ve-u) dd^c (v-u)\wedge (dd^c u)^{m-1}\wedge T'\\
&+\chi^{(m-1)} (\ve) \int_{V'} dd^c w_1 \wedge \cdots \wedge dd^c w_k.
\end{aligned}$$
It follows that 
$$\int_{V'} \chi \circ (v+\ve-u)T \wedge T' \le \int_{V'} w_1 \chi^{(m)} \circ (v+\ve-u) dd^c (v-u)\wedge (dd^c u)^{m-1}\wedge T'+
(\sum_{j=0}^{m-1} \chi^{(j)} (\ve))\int_{V'} T \wedge T'.$$
Since $u,v$ are bounded on $V'$, by letting $\ve \downarrow 0$ and applying Lebesgue dominated convergence's theorem 
(and taking into account the Chern-Levine-Nirenberg inequality (\ref{CLN}) and (\ref{chi})) we obtain
$$\int_{V'} \chi \circ (v-u) T \wedge T' \le 
\int_{V'} w_1 \chi^{(m)} \circ (v-u) dd^c (v-u) \wedge (dd^c u)^{m-1} \wedge T'+P_m (\chi) \int_{V'} T \wedge T'.$$
Since $u, v \in PSH(V)$ and $w_1<0$, the first term on the right hand side may be dominated as follows
$$\begin{aligned}
&\int_{V'} (-w_1) \chi^{(m)} \circ (v-u) dd^c (u-v) \wedge (dd^c u)^{m-1} \wedge T'\\
 &\le \int_{V'} (-w_1) \chi^{(m)} \circ (v-u)\big(\sum_{j=0}^{m-1} (dd^c u)^j \wedge (dd^c v)^{m-j-1}\big)
\wedge dd^c (u-v) \wedge T'\\
&=\int_{V'} (-w_1)\chi^{(m)} \circ (v-u)[(dd^c u)^m-(dd^c v)^m] \wedge T'.
\end{aligned}$$  
Putting all this together and rearranging we obtain (\ref{estimate}).

It remains to remove the restriction on smoothness of $\chi.$ Toward this end, we use Lemma 2.4 to get a sequence 
$\chi_j$ of $m-$increasing, $\mathcal C^\infty-$smooth functions such that $\chi_j$ and $\chi_j^{(l)}$ converges
locally uniformly to $\chi$ and $\chi^{(l)}$ on $[0, \infty)$ for $0 \le l \le m$. Then for each $j$, we have by (\ref{estimate})
$$\begin{aligned}
&\int_{V'} \chi_j \circ (v-u) T \wedge T' +
\int_{V'} (-w_1) \chi_j^{(m)} \circ (v-u)(dd^c v)^m \wedge T'\\
&\le \int_{V'} (-w_1) \chi_j^{(m)} \circ (v-u)(dd^c u)^m \wedge T'+P_m (\chi_j) \int_{V'} T \wedge T'.
\end{aligned}$$
By letting $j \to \infty$ and using Lebesgue dominated convergence's theorem we get (\ref{estimate}).

Finally, by letting $V' \uparrow V$ and using Lebesgue monotone convergence theorem in both sides of 
(\ref{estimate}) we complete the proof of the lemma.
 
\vskip0,3cm
\noindent
The final ingredient is the following equality of measures which is a modification of Proposition 4.2 in [BT2].
\vskip0,2cm
\noindent
{\bf Lemma 3.4.} {\it Let $1 \le m \le k$ and
$u, w_1, \cdots, w_{k-m} \in PSH(V) \cap L^\infty_{loc} (V), v \in PSH(V).$ Set
$T:= dd^c w_1 \wedge \cdots \wedge dd^c w_{k-m}.$ Then we have
$$(dd^c \max\{u,v\})^m \wedge T=(dd^c u)^m \wedge T \  \text{on}\  \{u>v\}.$$}
\noindent
{\it Proof.} We use some ideas in the proof of Theorem 4.1 in [KH].
Fix $a \in V$, it suffices to show that there is some open ball $\mathbb B\subset D$ around $a$ such that
$$(dd^c \max\{u,v\})^m \wedge T=(dd^c u)^m \wedge T \  \text{on}\  \{u>v\} \cap \mathbb B.$$
To see this, we first choose a small ball $\mathbb B \subset D$ around 
$a$ such that $u, w_1, \cdots, w_{k-m}$ are restrictions of \psh\ functions on $\mathbb B.$
By a standard smoothing process and shrinking $\mathbb B$, we may find
sequences of smooth \psh\ functions
$u_j, w_{1,j}, \cdots, w_{k-m,j}$ on $\mathbb B$ such that $u_j\downarrow u$ and 
$w_{1,j} \downarrow w_1, \cdots, w_{k-m, j}\downarrow w_{k-m}$ on $\mathbb B \cap V.$

Next we set 
$$T_j:= dd^c w_{1,j}\wedge \cdots dd^c w_{k-m,j}.$$
Since $\{u_j>v\} \cap V \cap \mathbb B$ is {\it open} in $V$, we have
$$(dd^c \max\{u_j, v\})^m \wedge T_j=(dd^c u_j)^m \wedge T_j \  \text{on}\ \{u_j>v\} \cap V \cap \mathbb B.$$
Since $\{u>v\} \cap \mathbb B \subset \{u_j>v\}$ we infer
$$(dd^c \max\{u_j,v\})^m \wedge T_j=(dd^c u_j)^m \wedge T_j \  \text{on}\ \{u>v\} \cap \mathbb B.$$
Set $u':=\max\{u-v, 0\}$. Then $u'$ is locally bounded and quasi-continuous on $V.$ 
Then we apply Proposition 2.2(c) to get the following weak convergences of measures on $V \cap \mathbb B$
$$\begin{aligned}
&u'(dd^c \max\{u_j, v\})^m \wedge T_j \to u' (dd^c \max\{u, v\})^m \wedge T,\\
&u'(dd^c u_j)^m \wedge T_j \to u' (dd^c u)^m \wedge T.
\end{aligned}$$
This implies that $u'\mu=0$ on $V \cap \mathbb B,$ where
$$\mu:= (dd^c \max\{u,v\})^m \wedge T-(dd^c u)^m \wedge T.$$
It follows, using Hahn's decomposition theorem for the real measure $\mu$ and the fact that $u'>0$ 
on $u>v$ (see Lemma 4.2 in [KH]), that $\mu=0$ on $\{u>v\}\cap \mathbb B.$ We are done.
\vskip0,3cm
\noindent
{\it Proof of Theorem 1.3.} First, we treat the case where $E=\emptyset$ and $u, v \in PSH(V) \cap L^\infty (V).$
For $\ve>0,$ we set
$v_\ve:=\max \{u, v-\ve\}.$ Then $v_\ve \in PSH(V) \cap L^\infty (V)$. Moreover, by the
the assumption (b) we see that $v_\ve=u$ on a \nhd\ of $\partial V$.
As in Lemma 3.3 we put
$$T:= dd^c w_1 \wedge \cdots \wedge dd^c w_m, T':=dd^c w_{m+1} \wedge \cdots \wedge dd^c w_k.$$
Then using Lemma 3.3 we get
$$
\int_V \chi \circ (v_\ve-u) T \wedge T' 
\le\int_V (-w_1)\chi^{(m)} \circ (v_\ve-u)\big [(dd^c u)^m-(dd^c v_\ve)^m\big]\wedge T'
+P_m (\chi) \int_V T \wedge T'.$$
By Lemma 3.4, we have 
$$\big [(dd^c u)^m-(dd^c v_\ve)^m\big]\wedge T'=0 \ \text{on}\
\{u>v-\ve\},\ \big [(dd^c v)^m-(dd^c v_\ve)^m\big]\wedge T'=0 \ \text{on}\ \{u<v-\ve\}.$$
Note also that 
$$v_\ve-u=v-\ve-u \ \text{on}\ u \le v-\ve.$$ 
From these facts we conclude that
$$
\begin{aligned}
&\int_{\{u<v-\ve\}} \chi \circ (v-\ve-u) T \wedge T'+
\int_{\{u<v-\ve\}} (-w_1) \chi^{(m)} \circ (v-\ve-u)(dd^c v)^m\wedge T'\\
&\le \int_{\{u \le v-\ve\}} (-w_1)\chi^{(m)} \circ (v-\ve-u)(dd^c u)^m\wedge T'+P_m (\chi) \int_V T \wedge T'.
\end{aligned}$$
Let $\theta_\ve$ and $\tilde \theta_\ve$ be the characteristic functions of $\{u<v-\ve\}$ and $\{u \le v-\ve\}$, respectively.
Then we see that $\theta_\ve \uparrow 1, \theta'_\ve \uparrow 1$ on $\{u<v\}$ as $\ve \downarrow 0.$
Hence, by applying Lebesgue's monotone convergence theorem and using the fact that $\chi$ is $m-$increasing 
we obtain 
$$\begin{aligned}
&\int_{\{u<v\}} \tilde \theta_\ve.(-w_1)\chi^{(m)} \circ (v-\ve-u) (dd^c v)^m\wedge T' \uparrow 
\int_{\{u<v\}} (-w_1) \chi^{(m)} \circ (v-u)(dd^c v)^m\wedge T',\\
&\int_{\{u<v\}} \theta_\ve.(-w_1)\chi^{(m)} \circ (v-\ve-u)(dd^c u)^m\wedge T' \uparrow
\int_{\{u<v\}} (-w_1)\chi^{(m)} \circ (v-u)(dd^c u)^m\wedge T',\\
&\int_{\{u<v\}} \theta_\ve.\chi \circ (v-\ve-u) T \wedge T' \uparrow
\int_{\{u<v\}} \chi \circ (v-u) T \wedge T'.
\end{aligned}$$
Putting all this together we obtain the desired conclusion.

For the general case we proceed as follows. Let $V_j \uparrow V$ be an increasing sequence of sub-domains in $V.$ 
Since $E$ is negligible we may find a function $\psi$ satisfying (\ref{negli}).
Fix $j \ge 1$. Set
$$v_j (z):= v(z)+\fr1{j}\psi(z)-\fr1{j}, \ \forall z \in V.$$
Then $v_j \uparrow v$ on $V.$
We claim that there exists $\al(j) \ge j$ such that 
$\varliminf\limits_{z \to \partial V_{\al(j)}} (u(z)-v_j(z)) \ge 0.$
Indeed, if this is false then we can find a sequence $z_k \to z^* \in \partial V$ such that for each $k$ we have
$$u(z_k)<v_j(z_k) =v(z_k)+\fr1{j} \psi(z_k)-\fr1{j}.$$
It implies, in view of the assumption (a), that
$$\fr1{j} \psi(z_k) \ge \inf_{k \ge 1} (u(z_k)-v(z_k)) +\fr1{j}>-\infty.$$
Hence, by (\ref{negli}) we must have $z^* \not\in E.$ On the other hand, by the condition (b) we get
$$0 \le \varliminf\limits_{k \to \infty} (u(z_k)-v(z_k)) \le -\fr{1}j,$$
which is clearly absurd. The claim follows. For simplicity of notation, we may assume that $\al(j)=j$ for every $j.$
Since $u, v_j \in PSH(V_j) \cap L^\infty (V_j)$ and satisfies 
$\varliminf\limits_{z \to \partial V_j} (u(z)-v_j(z)) \ge 0$ and since $(dd^c v_j)^m \wedge T' \ge (dd^c v)^m \wedge T'$ on $V$, by the result proved in the preceding case we get
$$
\begin{aligned}
&\int_{\{u<v_j\}} \chi \circ (v_j-u) T \wedge T'+
\int_{\{u<v_j\}} (-w_1) \chi^{(m)} \circ (v_j-u)(dd^c v)^m\wedge T'\\
&\int_{\{u<v_j\}} \chi \circ (v_j-u) T \wedge T'+
\int_{\{u<v_j\}} (-w_1) \chi^{(m)} \circ (v_j-u)(dd^c v_j)^m\wedge T'\\
&\le \int_{\{u<v_j\}} (-w_1)\chi^{(m)} \circ (v_j-u)(dd^c u)^m\wedge T'+P_m (\chi) \int_{V_j} dd^c w_1 \wedge \cdots \wedge dd^c w_k.
\end{aligned}$$
Observe that $\{z \in V_j: u(z)<v_j(z)\} \uparrow \{z \in V: u(z)<v(z)\}$ as $j \to \infty.$
So, by letting $j \to \infty$ and using Lebesgue monotone convergence theorem as in the previous case we obtain 
$$\begin{aligned}
&\int_{\{u<v\}} \chi \circ (v-u) T \wedge T'+
\int_{\{u<v\}} (-w_1) \chi^{(m)} \circ (v-u)(dd^c v)^m\wedge T'\\
&\le \int_{\{u<v\}} (-w_1)\chi^{(m)} \circ (v-u)(dd^c u)^m\wedge T'+P_m (\chi) \int_V T \wedge T'.
\end{aligned}$$
We have the desired result.
\vskip0,3cm
\noindent
{\it Proof of Corollary 1.4.} By applying Theorem 1.3 to $\chi(t)=t^m$ and 
$w_1=\cdots=w_k=\psi_r:=\fr{\Vert z\Vert^2}{r}-1$ with $r>\sup\{\Vert z\Vert^2: z \in V\}$ we obtain
$$\begin{aligned}
&C_{r,m}\int_{\{u<v\}} (v-u)^m \om^k+\int_{\{u<v\}} (-\psi_r) (dd^c v)^m \wedge \om^{k-m}\\
&\le \int_{\{u<v\}} (-\psi_r) (dd^c u)^m \wedge \om^{k-m},
\end{aligned}$$
where $C_{r,m}>0$ is a constant depends only on $r,m$.
By the assumption on $u,v$ we conclude that $$\int_{\{u<v\}} (v-u)^m \om^k=0.$$ This implies that $v \le u$ a.e. 
(with respect to $\om^k$) on $V_r$, the smooth locus of $V.$ 
Hence $v \le u$ entirely on $V_r.$
Now, we fix $a \in V_s$, we claim that there exists a one dimensional complex subvariety $\gamma \subset V$ such that 
$\gamma \cap V_s =\{a\}$. To see this, we first
make a change of coordinates to find a polydisc $\De$ in $\mathbb C^n$ that contains $a$ and a polydisc $\De'$ in 
$\mathbb C^k$ such that the projection map
$\pi: (z_1, \cdots, z_n) \mapsto (z_1, \cdots, z_k)$ expresses $V \cap \De$ as a 
branched cover of $\De'=\pi(\De)$ which is branched over a 
proper complex subvariety $H$ of $\De'$. 
Thus we can find a complex line $l \subset \mathbb C^k$ passing through $\pi(a)$ such that $l \cap H$ is discrete. 
Since $\pi(V_s \cap \De) \subset H$, we have $\gamma \cap V_s=\{a\}$, where $\gamma:=\pi^{-1} (U)$ and $U \subset l$ 
is a small \nhd\ of $\pi(a) \in l.$ This proves our claim. Next, we pick an irreducible branch $\gamma' \subset \gamma$ that contains $a.$
Then, by normalization we can find a connected Riemann surface $\gamma^*$ and holomorphic mapping 
$f: \gamma^* \to \gamma$ which is surjective. 
Set $u':=u \circ f|_{\gamma'}, v':=v \circ f|_{\gamma'}.$ Since $u \le v$ on $V_r,$ by the choice of $\gamma$, we infer that $v' \le u'$ on $\gamma^*$ except for the finite set $f^{-1} (a)$. Hence, this inequality holds true entirely on $\gamma^*$ since $u', v'$ are subharmonic there. It follows that $v(a) \le u(a).$
The proof is complete.
\vskip0,4cm
\noindent
{\it Proof of Corollary 1.5.} Define inductively on $(0, \infty)$ the following functions 
$$\chi_0 (t):=\chi (t), \chi_{j+1} (t)=\int_0^t \chi_j (x)dx, j \ge 0.$$
We can check that $\chi_k^{(j)}=\chi_{k-j}$ for every $0 \le j \le k.$ 
In particular $\chi_k^{(k)}=\chi_0$ on $(0, \infty)$ and $\chi_k$ is $k-$increasing. Furthermore, $P_k (\chi_k)=0$. 
Now we apply Theorem 1.3 to the "weight" function $\chi_k$ and
$w_1=\cdots=w_k=\psi_r:=\fr{\Vert z\Vert^2}{r}-1$ with $r>\sup\{\Vert z\Vert^2: z \in V\}$ to obtain
$$\int_{\{u<v\}} (-\psi_r) \chi \circ (v-u)(dd^c v)^k \le \int_{\{u<v\}} (-\psi_r) \chi \circ (v-u)(dd^c u)^k.$$
By letting $r \to \infty$ we arrived at the desired conclusion.
\vskip0,4cm
\noindent
{\it Proof of Corollary 1.6.} First, we let $\chi_k$ be the function constructed in the proof of Corollary 1.5.
Next, fix $\de>0$. Set
$$A_j:= \{z \in V: u(z)>u_j(z)+\de\}, B_j:= \{z \in V: u_j(z)>u(z)+\de\}.$$
We claim that $C(A_j, V) \to 0$ as $j \to \infty.$ Assume otherwise, then, by switching to a subsequence, we may find a sequence $\{\psi_j\} \subset PSH(V), -1<\psi_j<0$ and $\la>0$ such that
$$\int_{A_j} (dd^c \psi_j)^k \ge \la, \ \forall j.$$
Fix $j \ge 1.$ In view of the assumption (a), we may apply Theorem 1.3 to $u, u_j, w_1=\cdots=w_k=\psi_j$ and 
$\chi_k$ to obtain
$$\begin{aligned}
\int_{\{u_j<u\}} \chi_k \circ (u-u_j) (dd^c \psi_j)^k  &\leq
\int_{\{u_j<u\}} (-\psi_j)\chi\circ (u-u_j) d\mu_j\\
&\le \int_{\{u_j<u\}} \chi\circ (u-u_j) d\vert \mu_j\vert.
\end{aligned}$$ 
It implies, using the condition (b), that
$$\lim_{j \to \infty} \int_{A_j} \chi_k \circ (u-u_j) (dd^c \psi_j)^k =0.$$
On the other hand, for each $j \ge 1$ we have
$$\int_{A_j} \chi_k \circ (u-u_j) (dd^c \psi_j)^k \ge \chi_k (\de) \int_{A_j} (dd^c \psi_j)^k \ge 
\la \chi_k (\de)>0.$$
We arrived at a contradiction. Hence $\lim\limits_{j \to \infty} C(A_j,V)=0.$ 
By exchanging the role of $u$ and $u_j$ and repeating the same reasoning we also obtain
$\lim\limits_{j \to \infty} C(B_j, V)=0.$ 
The proof is thereby completed.
\vskip1cm
\cen{\bf REFERENCES}
\vskip0,5cm
\noindent
[Be] E. Bedford, {\it The operator $(dd^c)^n$ on complex spaces}, S\'eminaire Lelong-Skoda, Springer Lecture Notes {\bf 919} (1981), 294-323.

\noindent
[BT1] E. Bedford and A. Taylor, {\it A new capacity for \psh\ functions,} Acta. Math., {\bf 149} (1982), 1-40.
 
\noindent
[BT2] E. Bedford and A. Taylor, {\it Fine topology, Shilov boundary, and $(dd^c)^n,$} Journal of Functional Analysis {\bf 72},
(1987), 225-251.

\noindent
[BT3] E. Bedford and A. Taylor, {\it Uniqueness for the complex Monge-Amp\`ere equation for functions with logarithmic growth}, Indiana Univ. Math. J., {\bf 38} (1989), 455-469.

\noindent
[Ce] U. Cegrell, {\it The general definition of the complex Monge-Amp\`ere operator}, Ann.
Inst. Fourier (Grenoble) {\bf 54} (2004), 159-179.

\noindent
[FN] J. E. Fornaess and R. Narasimhan, {\it The Levi problem on complex spaces with singularities},
Math. Ann. {\bf 248} (1980) 47-72.

\noindent
[GH] P. Griffiths and J. Harris, {\it Principles of Algebraic Geometry}, Wiley Classics Library, John
Wiley $\&$ Sons, Inc., New York, 1994. 

\noindent
[KH] N.V. Khue and P.H. Hiep, {\it A comparison principle for the complex Monge-Amp\`ere operator in Cegrell's classes
and applications}, Trans. Amer. Math. Soc., {\bf 361} (2009), 5539-5554.

\noindent
[Xi1] Y. Xing, {\it Continuity of the complex Monge-Amp\`ere operator}, Proc. Amer. Math. Soc., {\bf 124} (1996), 457-467.

\noindent
[Xi2] Y. Xing, {\it A strong comparison principle for \psh\ functions with finite pluricomplex energy}, Michigan Math.
J., {\bf 56} (2008), 563-581.  
\end{document}